\documentclass[11pt]{amsart}
\usepackage{amssymb}
\usepackage{amsmath}
\usepackage[active]{srcltx}
\usepackage{t1enc}
\usepackage[latin2]{inputenc}
\usepackage{verbatim}
\usepackage{amsmath,amsfonts,amssymb,amsthm}
\usepackage[mathcal]{eucal}
\usepackage{enumerate}
\usepackage[centertags]{amsmath}
\usepackage{graphics}

\newtheorem{theorem}{Theorem}

\newtheorem{lemma}{Lemma}
\newtheorem{remark}{Remark}

\newtheorem{corollary}{Corollary}

\begin{document}
\author{Nika Areshidze, Davit Baramidze, Lars-Erik Persson, George Tephnadze}
\title[Fejér means]{Weighted maximal operators of Fejér means of Walsh-Fourier series in the martingale Hardy space $H_{1/2}$}

\address{N. Areshidze, Department of Mathematics, Faculty of Exact and Natural Sciences, Tbilisi State University, Chavchavadze str. 1, Tbilisi 0128, Georgia}
\email{nika.areshidze15@gmail.com }
\address{D. Baramidze, The University of Georgia, School of science and
	technology, 77a Merab Kostava St, Tbilisi 0128, Georgia and Department of
	Computer Science and Computational Engineering, UiT - The Arctic University
	of Norway, P.O. Box 385, N-8505, Narvik, Norway.}
\email{davit.baramidze@ug.edu.ge }
\address{Lars-Erik Persson, UiT The Arctic University of Norway, P.O. Box 385, N-8505, Narvik, Norway and Department of Mathematics and Computer Science, Karlstad University, 65188 Karlstad, Sweden.}
\email{larserik6pers@gmail.com}
\address{G. Tephnadze, The University of Georgia, School of science and
	technology, 77a Merab Kostava St, Tbilisi 0128, Georgia.}
\email{g.tephnadze@ug.edu.ge}

\thanks{The research was supported by Shota Rustaveli National Science Foundation grant no. PHDF-21-1702.}
\date{}
\maketitle

\begin{abstract}
In this paper we derive the restricted weighted maximal operator, defined by
${\sup }_{k\in \mathbb{N}}\left(\left\vert \sigma _{k}F\right\vert/A^2_k\right)$ 
of Fejér means of Walsh-Fourier series and prove that the it is bounded  from the martingale Hardy space $H_{1/2}(G)$ to the Lebesgue space $L_{1/2}(G).$ The sharpness of this result is also proved. As a consequence we obtain some new and and well-know results.
\end{abstract}

\date{}



\section{INTRODUCTION}

\bigskip In the one-dimensional case the weak (1,1)-type inequality for the
maximal operator of Fejér means with respect to Walsh system, defined by
\begin{equation*}
\sigma ^{\ast }f:=\sup_{n\in \mathbb{N}}\left\vert \sigma _{n}f\right\vert
\end{equation*}%
was investigated in Schipp \cite{Sc} and Pál, Simon \cite{PS} (see also \cite{BNPT}, \cite{NTT} and \cite{PTW2}). Fujii \cite{Fu} and Simon \cite{Si2} verified
that $\sigma ^{\ast }$ is bounded from $H_{1}(G)$ to $L_{1}(G)$. Weisz \cite{We2}
generalized this result and proved boundedness of $\sigma ^{\ast }$ from the
martingale space $H_{p}$ to the Lebesgue space $L_{p}$ for $p>1/2$. Simon 
\cite{Si1} gave a counterexample, which shows that boundedness does not hold
for $0<p<1/2.$ A counterexample for $p=1/2$ was given by Goginava 
\cite{GoAMH}. Moreover, \cite{Goginava}  (see also \cite{tep1}) he  proved that there
exists a martingale $F\in H_{p}(G)$ $\left( 0<p\leq 1/2\right) ,$ such that 
$$
\sup_{n\in \mathbb{N}}\left\Vert \sigma _{n}F\right\Vert _{p}=+\infty .
$$
Weisz \cite{We4} proved that
the maximal operator $\sigma ^{\ast }$  of the Fejér means is bounded from
the Hardy space $H_{1/2}$ to the space $weak-L_{1/2}$. 
Weisz \cite{We3} (see also \cite{We1}) also proved that for any $F\in H_p(G)$ $(p>0),$  the maximal operator 
$$
\underset{n\in \mathbb{N}}{\sup }\left\vert \sigma
_{2^n}F\right\vert
$$
is bounded from
the Hardy space $H_{p}(G)$ to the space $L_{p}(G)$, but (for details see \cite{PTWbook}) it is not bounded from the hardy space  $H_{p}(G)$ to the space $H_{p}(G).$

To study convergence of subsequences of  Fejér means and their restricted  maximal operators on the martingale Hardy spaces $H_p(G)$ for $0<p\leq 1/2,$ the central role is played by the fact that any natural number $n\in \mathbb{N}$ can be uniquely expression as
$$ 
n=\sum_{k=0}^{\infty }n_{j}2^{j},  \ \ n_{j}\in Z_{2} ~(j\in \mathbb{N}), $$ 
where only a finite numbers of $n_{j}$ differ from zero
and their important characters  $\left[ n\right],$ $\left\vert n\right\vert,$ $\rho\left( n\right)$  and $V(n)$ are defined by
\begin{equation*}
\left[ n\right] :=\min \{j\in \mathbb{N},n_{j}\neq 0\}, \ \ 
\left\vert n\right\vert :=\max \{j\in \mathbb{N},n_{j}\neq 0\}, \ \ 
\rho\left( n\right) =\left\vert n\right\vert -\left[ n\right] 
\end{equation*}
and
\begin{equation*}
V\left( n\right): =n_{0}+\overset{\infty }{\underset{k=1}{\sum }}\left|
n_{k}-n_{k-1}\right|, \text{ \ for \
	all \ \ }n\in \mathbb{N}.
\end{equation*}

Persson and Tephnadze \cite{PT} generalized this result and proved that if  $0<p\leq 1/2$ and $\left\{ n_{k}:k\geq 0\right\} $ is a sequence of positive numbers, such that
$$\sup_{k\in \mathbb{N}}\rho \left( n_{k}\right) \leq c<\infty,$$
then the maximal operator $\widetilde{\sigma }^{\ast ,\nabla },$ defined by
\begin{equation}\label{max}
\widetilde{\sigma }^{\ast ,\nabla }F:=\underset{k\in \mathbb{N}}{\sup }
\left\vert \sigma _{n_{k}}F\right\vert
\end{equation}%
	is bounded from the Hardy space $H_{p}(G)$ to the Lebesgue space $L_{p}(G).$
	Moreover, if $0<p<1/2$ and $\left\{ n_{k}:k\geq 0\right\} $ is a sequence of positive numbers, such that
$$
\sup_{k\in \mathbb{N}}\rho \left( n_{k}\right) =\infty , 
$$
then there exists a martingale $F\in H_{p}(G)$ such that 
$$
\underset{k\in \mathbb{N}}{\sup }\left\Vert \sigma _{n_{k}}F\right\Vert_{p}=\infty .
$$
In \cite{tep8} it was proved that if $F\in H_{1/2}(G),$ then there exists an absolute
	constant $c,$ such that%
	\begin{equation*}
	\left\Vert \sigma _{n}F\right\Vert _{H_{1/2}}\leq cV^{2}\left( n\right)
	\left\Vert F\right\Vert _{H_{1/2}}.
	\end{equation*}	
Moreover, if  $\left\{ n_{k}:k\geq 0\right\} $ is subsequence of positive integers $\mathbb{N}_{+},$ such that 
$$\sup_{k\in \mathbb{N}}V\left(n_{k}\right) =\infty $$	
and  $\{\varphi_n\}$ is any nondecreasing, nonnegative function, satisfying conditions $\Phi_n \uparrow \infty $ and
\begin{equation*}
\overline{\underset{k\rightarrow \infty }{\lim }}\frac{V^{2}\left(n_{k}\right) }{\Phi_{n_{k} }}=\infty , 
\end{equation*}
then there exists a martingale $F\in H_{1/2}(G),$ such that
	\begin{equation*}
	\underset{k\in \mathbb{N}}{\sup }\left\Vert \frac{\sigma _{n_{k}}F}{\Phi_{n_{k}} }\right\Vert _{1/2}=\infty .
	\end{equation*}
It follows that if $f\in H_{1/2}(G)$ and $\left\{ n_{k}:k\geq 0\right\} $ is any sequence of positive numbers, then $\sigma _{n_{k}}f$ are bounded	from the Hardy space $H_{1/2}(G)$ to the space $H_{1/2}(G)$ if and only if, for some $c,$
$$\sup_{k\in \mathbb{N}}V\left(n_{k}\right)<c<\infty.$$

In \cite{tep2} it was proved that the weighted  maximal operator
\begin{equation*}
\overset{\sim }{\sigma }^{*}f:=\sup_{n\in \mathbb{N}}\frac{\left| \sigma _{n}f\right|
}{\log ^{2}\left( n+1\right) }
\end{equation*}
is bounded from the Hardy space $H_{1/2}\left( G\right) $ to the space $
L_{1/2}\left( G\right) .$
Moreover, it was also proved that the rate of denominator $\{\log^2(n+1)\}$ can not be improved.

Baramidze and Tephnadze \cite{BaTe1} proved that if  $\left\{ n_{k}:k\geq
0\right\} $ be a sequence of positive numbers, such that
$$
\sup_{k\in \mathbb{N}}A_{ \vert n_{k}\vert } \leq c<\infty,
$$
then the maximal operator \eqref{max} is bounded from the Hardy space $H_{1/2}(G)$ to the space $L_{1/2}(G).$ Moreover, if $\left\{ n_{k}:k\geq 0\right\} $ is a sequence of positive numbers, such that
$$
\sup_{k\in \mathbb{N}}A_{\vert n_{k}\vert} =\infty , 
$$
then there exists a martingale $F\in H_{p}(G)$ such that 
$$
\underset{k\in \mathbb{N}}{\sup }\left\Vert \sigma _{n_{k}}F\right\Vert_{1/2}=\infty .
$$

In this paper we investigate the restricted weighted maximal operator, defined by 
$${\sup }_{k\in \mathbb{N}}\left(\left\vert \sigma _{k}F\right\vert/A^2_k\right)$$
of Fejér means of Walsh-Fourier series and prove that the it is bounded  from the martingale Hardy space $H_{1/2}(G)$ to the Lebesgue space $L_{1/2}(G).$ As a consequence we obtain some new and and well-know results.

This paper is organized as follows: In order not to disturb our discussions later on some
definitions and notations are presented in Section 2. The main result and some of its
consequences can be found in Section 3. For the proof of the main result we need some auxiliary statements. These results are
presented in Section 4. The detailed proofs are given in Section 5.

\section{Definitions and Notations}

Let $\mathbb{N}_{+}$ denote the set of the positive integers, $\mathbb{N}:=
\mathbb{N}_{+}\cup \{0\}.$ Denote by $Z_{2}$ the discrete cyclic group of
order 2, that is $Z_{2}:=\{0,1\},$ where the group operation is the modulo 2
addition and every subset is open. The Haar measure on $Z_{2}$ is given so
that the measure of a singleton is 1/2.

Define the group $G$ as the complete direct product of the group $Z_{2},$
with the product of the discrete topologies of $Z_{2}$. The elements of $G$
are represented by sequences $x:=(x_{0},x_{1},...,x_{j},...),$ where $%
x_{k}=0\vee 1.$

It is easy to give a base for the neighborhood of $x\in G$ 
\begin{equation*}
I_{0}\left( x\right) :=G,\text{ \ }I_{n}(x):=\{y\in
G:y_{0}=x_{0},...,y_{n-1}=x_{n-1}\}\text{ }(n\in \mathbb{N}).
\end{equation*}

Denote $I_{n}:=I_{n}\left( 0\right) ,$ $\overline{I_{n}}:=G$ $\backslash $ $%
I_{n}$ and $e_{n}:=\left( 0,...,0,x_{n}=1,0,...\right) \in G,$ for $n\in 
\mathbb{N}$. Then it is easy to prove that 
\begin{equation}\label{1}
\overline{I_{M}}=\underset{i=0}{\bigcup\limits^{M-1}}I_{i}\backslash
I_{i+1}=\left( \overset{M-2}{\underset{k=0}{\bigcup }}\overset{M-1}{\underset%
{l=k+1}{\bigcup }}I_{l+1}\left( e_{k}+e_{l}\right) \right) \bigcup \left( 
\underset{k=0}{\bigcup\limits^{M-1}}I_{M}\left( e_{k}\right) \right).
\end{equation}

If $n\in \mathbb{N},$ then every $n$ can be uniquely expressed as $%
n=\sum_{j=0}^{\infty }n_{j}2^{j},$ where $n_{j}\in Z_{2}$ $~(j\in \mathbb{N})
$ and only a finite numbers of $n_{j}$ differ from zero.
Every $n\in \mathbb{N}$ can be also represented as $n=%
\sum_{i=1}^{r}2^{n_{i}},n_{1}>n_{2}>...n_{r}\geq 0.$ For such a representation
of $n\in \mathbb{N},$ we denote numbers 
\begin{equation*}
n^{\left( i\right) }=2^{n_{i+1}}+...+2^{n_{r}},i=1,...,r.
\end{equation*}
Let 
$\begin{matrix}
2^{s}\le {{n}_{{{s}_{1}}}}\le {{n}_{{{s}_{2}}}}\le ...\le {{n}_{{{s}_{r}}}}\le {2^{s+1}}, \ s\in \mathbb{N}.  \\
\end{matrix}$
For such ${{n}_{{{s}_{j}}}},$ which can be written as
${{n}_{{{s}_{j}}}}=\sum\limits_{i=1}^{{{r}_{{{s}_{j}}}}}{\sum\limits_{k=l_{i}^{{{s}_{j}}}}^{t_{i}^{{{s}_{j}}}}{2^k}},$
where   
$0\le l_{1}^{{{s}_{j}}}\le t_{1}^{{{s}_{j}}}\le l_{2}^{{{s}_{j}}}-2<l_{2}^{{{s}_{j}}}\le t_{2}^{{{s}_{j}}}\le ...\le l_{{{r}_{j}}}^{{{s}_{j}}}-2<l_{{{r}_{s_j}}}^{{{s}_{j}}}\le t_{{{r}_{s_j}}}^{{{s}_{j}}},$
we denote
\begin{eqnarray*}
{{A}_{s}}&:=&\bigcup\limits_{j=1}^{r}{\left\{ l_{1}^{{{s}_{j}}},t_{1}^{{{s}_{j}}},l_{2}^{{{s}_{j}}},t_{2}^{{{s}_{j}}},...,l_{{{r}_{s_j}}}^{{{s}_{j}}}, t_{{{r}_{s_j}}}^{{{s}_{j}}} \right\}}\\ \notag
&=&\bigcup\limits_{j=1}^{r}{\left\{ l_{1}^{{{s}_{j}}},l_{2}^{{{s}_{j}}},...,l_{{{r}_{s_j}}}^{{{s}_{j}}} \right\}}
\bigcup\limits_{j=1}^{r}{\left\{ t_{1}^{{{s}_{j}}},t_{2}^{{{s}_{j}}},...,t_{{{r}_{s_j}}}^{{{s}_{j}}} \right\}}\\ \notag
&=&{\left\{ l_{1}^{s},l_{2}^{s},...,l_{r^1_{s}}^{s} \right\}}
\bigcup{\left\{ t_{1}^{s},t_{2}^{s},...,t_{{{r}_{s}^2}}^{s} \right\}}
={\left\{ u_{1}^{s},u_{2}^{s},...,u_{r^3_{s}}^{s} \right\}},
\end{eqnarray*}
where $ u_{1}^{s}<u_{2}^{s}<...<u_{r^3_{s}}^{s}.$
We note that $
t_{{{r}_{s_j}}}^{{{s}_{j}}}=s\in {{A}_{s}}, \ \ \text{ for } \ \  j=1,2,...,r.  
$

Let us denote the cardinality of the set $A_s$ by $\vert A_s\vert$, that is
$$card(A_s):=\vert A_s\vert.$$ 
It is evident that
$
\vert A_s\vert=r^3_s\leq r_s^1+r_s^2.
$
Moreover,
$r^s_3=card(A_s)<\infty$
if and only if
$r_s^1<\infty \ \  \text{and} \ \ r_s^2<\infty.$
We note that if $\vert A_s\vert<\infty,$ then each ${{n}_{{{s}_{j}}}}$ has bounded variation
\begin{equation}\label{cond1x}
V(n_{s_j})<c<\infty, \ \ \ \text{for each}  \ \ \ j=1,2,\ldots,r.
\end{equation}
and  \ \
$r_s^1<\infty, \  r_s^2<\infty \   \text{and} \  r_s^3<\infty.$
Therefore, if we consider blocks (intervals)
\begin{eqnarray}\label{000}
&& [u_1^s,u_1^s],  [u_1^s,u_2^s], ..., [u_1^s,u_{r_s^3}^s], [u_2^s,u_2^s],  [u_2^s,u_3^s], ..., [u_2^s,u_{r_s^3}^s]\ldots
[u_{r_s^3}^s,u_{r_s^3}^s],
\end{eqnarray}
then it is easy to see that it contains  $\left({{r_s^3}{(r_s^3+1)}}\right)/{2}$ different blocks. Therefore, the dyadic representation of different natural numbers, which contains blocks from \eqref{000},
can be at most $2^{\frac{{r_s^3}{(r_s^3+1)}}{2}}-1,$ which is finite number and the set
$\{{{n}_{{{s}_{1}}}},{{n}_{{{s}_{2}}}}, ...{{n}_{{{s}_{r}}}}\}$
is finite for all $s\in \mathbb{N}_+,$ from which it follows that 
\begin{equation}\label{cond2x}
s_r<\infty, \ \ \ \text{for all}\ \ \ s\in\mathbb{N}.
\end{equation}
Summing up, we can conclude that
$\sup_{s\in\mathbb{N}}\vert A_s\vert<\infty$
if and only if the set  $\{{{n}_{{{s}_{1}}}},{{n}_{{{s}_{2}}}}, ...,{{n}_{{{s}_{r}}}}\}$ is finite for all $s\in \mathbb{N}_+$ and each ${{n}_{{{s}_{j}}}}$ has bounded variation, that is, conditions \eqref{cond1x} and \eqref{cond2x} are fulfilled.

The norms (or quasi-norm) of the spaces $L_{p}$ and $\text{weak}-L_{p },\left( 0<p<\infty \right) $ are respectively defined by 
\begin{equation*}
\left\Vert f\right\Vert _{p}^{p}:=\int_{G}\left\vert f\right\vert ^{p}d\mu
\ \ \ \ \text{and} \ \ \ \ \left\Vert f\right\Vert _{\text{weak}-L_{p}}^{p}:=\sup_{\lambda
>0}\lambda ^{p}\mu \left( f>\lambda \right).
\end{equation*}

The $k$-th Rademacher function is defined by 
\begin{equation*}
r_{k}\left( x\right) :=\left( -1\right) ^{x_{k}}\text{\qquad }\left(
x\in G,\text{ }k\in \mathbb{N}\right) .
\end{equation*}

Now, define the Walsh system $w:=(w_{n}:n\in \mathbb{N})$ on $G$ as: 
\begin{equation*}
w_{n}(x):=\overset{\infty }{\underset{k=0}{\Pi }}r_{k}^{n_{k}}\left(
x\right) =r_{\left\vert n\right\vert }\left( x\right) \left( -1\right) ^{%
\underset{k=0}{\overset{\left\vert n\right\vert -1}{\sum }}n_{k}x_{k}}\text{%
\qquad }\left( n\in \mathbb{N}\right) .
\end{equation*}

The Walsh system is orthonormal and complete in $L_{2}\left( G\right) $ (see 
\cite{sws}).

If $f\in L_{1}\left( G\right) ,$ then we can define the Fourier coefficients,
partial sums of Fourier series, Fejér means, Dirichlet and Fejér kernels in
the usual manner: 
\begin{eqnarray*}
\widehat{f}\left( n\right) &:=&\int_{G}fw_{n}d\mu ,\,\,\left( n\in \mathbb{N}%
\right),\\
S_{n}f&:=&\sum_{k=0}^{n-1}\widehat{f}\left( k\right) w_{k},%
\text{\ }\left( n\in \mathbb{N}_{+},S_{0}f:=0\right) ,
\\
\sigma _{n}f&:=&\frac{1}{n}\sum_{k=1}^{n}S_{k}f,\text{ \ \ } \\
D_{n}&:=&\sum_{k=0}^{n-1}w_{k\text{ }},\text{ \ }\\
K_{n}&:=&\frac{1}{n}\overset{n}{%
\underset{k=1}{\sum }}D_{k}\text{ },\text{ \ }\left(n\in \mathbb{N}
_{+}\right) .
\end{eqnarray*}

Recall that (see \cite{G-E-S}  and \cite{sws})
\begin{equation}\label{1dn}
D_{2^{n}}\left( x\right) =\left\{ 
\begin{array}{ll}
2^{n} & \,\text{if\thinspace \thinspace \thinspace }x\in I_{n} \\ 
0\, & \ \,\text{if}\,\,x\notin I_{n}.%
\end{array}%
\right.  
\end{equation}

Let 
$n=\sum_{i=1}^{r}2^{n_{i}}, \ \ \ n_{1}>n_{2}>...>n_{r}\geq 0.$
Then (see \cite{gat1} and \cite{sws}) 
\begin{equation*}
nK_{n}=\sum_{A=1}^{r}\left( \underset{j=1}{\overset{A-1}{\prod }}%
w_{2^{n_{j}}}\right) \left( 2^{n_{A}}K_{2^{n_{A}}}+n^{\left( A\right)
}D_{2^{n_{A}}}\right).  
\end{equation*}

The $\sigma $-algebra, generated by the intervals $\left\{ I_{n}\left(
x\right) :x\in G\right\} $ will be denoted by $\zeta _{n}$ $\left( n\in 
\mathbb{N}\right) .$ Denote by $F=\left( F_{n},n\in \mathbb{N}\right) $ a
martingale with respect to $\zeta _{n}$ $\left( n\in \mathbb{N}\right) $
(for details see e.g. \cite{We1}). The maximal function $F^{\ast }$ of a martingale $F$
is defined by 
$
F^{\ast }:=\sup_{n\in \mathbb{N}}\left\vert F_{n}\right\vert .
$

In the case $f\in L_{1}\left( G\right) ,$ the maximal functions $f^{\ast }$ are also given by 
\begin{equation*}
f^{\ast }\left( x\right) :=\sup\limits_{n\in \mathbb{N}}\left( \frac{1}{\mu
\left( I_{n}\left( x\right) \right) }\left\vert \int_{I_{n}\left( x\right)
}f\left( u\right) d\mu \left( u\right) \right\vert \right) .
\end{equation*}

For $0<p<\infty ,$ the Hardy martingale spaces $H_{p}\left( G\right) $
consist of all martingales, for which 
\begin{equation*}
\left\Vert F\right\Vert _{H_{p}}:=\left\Vert F^{\ast }\right\Vert
_{p}<\infty .
\end{equation*}

A bounded measurable function $a$ is a $p$-atom, if there exists an interval $I,$ such that
\begin{equation*}
\text{ supp}\left( a\right) \subset I, \ \ \ \int_{I}ad\mu =0,\text{ \ }\left\Vert a\right\Vert _{\infty }\leq \mu \left(
I\right) ^{-1/p}.
\end{equation*}

It is easy to check that for every martingale $F=\left( F_{n},n\in \mathbb{N}%
\right) $ and every $k\in \mathbb{N}$ the limit
\begin{equation*}
\widehat{F}\left( k\right) :=\lim_{n\rightarrow \infty }\int_{G}F_{n}\left(
x\right) w_{k}\left( x\right) d\mu \left( x\right)
\end{equation*}%
exists and it is called the $k$-th Walsh-Fourier coefficients of $F.$

The Walsh-Fourier coefficients of $f\in L_{1}\left( G\right) $ are the same
as those of the martingale $\left( S_{2^{n}}f, n\in \mathbb{N}%
\right) $ obtained from $f$.

\section{The Main Result and its Consequences}

Our main result reads:
\begin{theorem}\label{theorem1}
a) Let $f\in {{H}_{1/2}}\left(G \right)$ and $\left\{ n_{k}:k\geq 0\right\} $ is a sequence of positive numbers. Then the weighted maximal operator $\widetilde{\sigma }^{\ast ,\nabla },$ defined by
\begin{equation*}
\widetilde{\sigma }^{\ast ,\nabla }F:=\underset{k\in \mathbb{N}}{\sup }
\frac{\left\vert \sigma _{n_k}F\right\vert}{A^2_{\vert n_ k\vert }},
\end{equation*}
is bounded from the Hardy space ${{H}_{1/2}}$ to the Lebesgue space ${{L}_{1/2}}$.

b) (Sharpness) Let  
\begin{equation}\label{cond2}
{{\sup }_{k\in \mathbb{N}}}\vert{{A}_{n_k}}\vert=\infty
\end{equation}
and 
$\{\varphi_n\}$ is a nondecreasing sequence
satisfying the condition
\begin{equation*}
\overline{\lim_{k\rightarrow \infty }}\frac{A_{\vert n_k\vert}^{2} }{
	\varphi_{\vert n_k\vert}}=\infty .
\end{equation*}
Then there exists a martingale $f\in {{H}_{1/2}}\left(G \right),$ such that the maximal operator
\begin{equation*}
\underset{k\in \mathbb{N}}{\sup }
\frac{\left\vert \sigma _{n_k}F\right\vert}{\varphi_{\vert n_k\vert}}
\end{equation*}
  is not bounded from the Hardy space ${{H}_{1/2}}(G)$ to the Lebesgue space ${{L}_{1/2}}(G).$
\end{theorem}

\begin{corollary}
	The  maximal operator $\overset{\sim }{\sigma }^{*},$ defined by
\begin{equation*}
\overset{\sim }{\sigma }^{*}f:=\sup_{n\in \mathbb{N}}\frac{\left| \sigma _{n}f\right|
}{\log ^{2}\left( n+1\right) }
\end{equation*}
is bounded from the Hardy space $H_{1/2}\left( G\right) $ to the space $
L_{1/2}\left( G\right) .$

b) Let $\{\varphi_n\}$ is a nondecreasing sequence satisfying the condition
	\begin{equation*}
	\overline{\lim_{n\rightarrow \infty }}\frac{\log ^{2}\left( n+1\right) }{
		\varphi_n}=+\infty .
	\end{equation*}
	Then there exists a martingale $f\in H_{1/2},$ such that the maximal operator
	\begin{equation*}
\sup_{n\in \mathbb{N}}\frac{\left\vert\sigma _{n}f\right\vert}{\varphi_n }
	\end{equation*}
is not bounded from the Hardy space ${{H}_{1/2}}(G)$ to the Lebesgue space ${{L}_{1/2}}(G).$
\end{corollary}

In order to be able to compare with some other results in the literature (see Remark \ref{theorem4}) we also state the following:
\begin{corollary}\label{theorem2}
Let $f\in {{H}_{1/2}}\left( G \right)$. Then the restricted maximal operators $\widetilde{\sigma }^{\ast ,\nabla }_i, \ i=1,2,3,$ defined by 
\begin{equation}\label{maxoperator00}
\widetilde{\sigma }^{\ast ,\nabla }_1F:=\underset{k\in \mathbb{N}}{\sup }
\left\vert \sigma _{2^k}F\right\vert,
\end{equation}
\begin{equation}\label{maxoperator0}
\widetilde{\sigma }^{\ast ,\nabla }_2F:=\underset{k\in \mathbb{N}}{\sup }
\left\vert \sigma _{2^k+1}F\right\vert,
\end{equation}
	\begin{equation}\label{maxoperator1}
	\widetilde{\sigma }^{\ast ,\nabla }_3F:=\underset{k\in \mathbb{N}}{\sup }
	\left\vert \sigma _{2^k+2^{[k/2]}}F\right\vert,
	\end{equation}
where $[n]$ denotes integer part of $n$, all	are bounded  from the Hardy space ${{H}_{1/2}}(G)$ to the Lebesgue space ${{L}_{1/2}}(G)$.
\end{corollary}
\begin{remark}\label{theorem4}
In \cite{PT} it was proved that if  $0<p<1/2,$ then the restricted maximal operators $\widetilde{\sigma }^{\ast ,\nabla }_2$ and $\widetilde{\sigma }^{\ast ,\nabla }_3$ defined by \eqref{maxoperator0} and \eqref{maxoperator1},
	are not bounded  from the Hardy space ${{H}_{p}}(G)$ to the Lebesgue space $weak-{{L}_{p}}(G)$.
	
On the other hand, in \cite{PT} it was proved that if $0<p\leq 1/2$, then the restricted maximal operator   $\widetilde{\sigma }^{\ast ,\nabla }_1$ defined by \eqref{maxoperator00} is bounded  from the Hardy space ${{H}_{p}}(G)$ to the Lebesgue space $L_p(G).$ 
\end{remark}

\begin{corollary}\label{theorem5}
Let $f\in {{H}_{1/2}}\left( G \right)$ and 
$\left\{ A^s_{k}:0\leq k\leq s-1\right\} $ be a sequence of positive numbers, defined by
$$\alpha^s_{0}=2^s+2^0, \ \ \ \alpha^s_{1}=2^s+2^0+2^1, ..., \alpha^s_{s-1}=2^s+2^0+2^1+...+2^{s-1}, \ \ \ s\in\mathbb{N}.$$
Then the restricted maximal operator $\widetilde{\sigma }^{\ast ,\nabla }_4,$ defined by 
\begin{equation*}
\widetilde{\sigma }^{\ast ,\nabla }_4F:=\underset{s\in \mathbb{N}}{\sup }\underset{0\leq k\leq s-1}{\sup }\frac{
\left\vert \sigma _{\alpha^s_k}F\right\vert}{\log^2{\alpha^s_k}},
\end{equation*}
is bounded  from the Hardy space ${{H}_{1/2}}$ to the Lebesgue space ${{L}_{1/2}}.$
\end{corollary}
\begin{corollary}\label{theorem6}
	Let $f\in {{H}_{1/2}}\left( G \right)$ and 
	$\left\{ A^s_{k}:0\leq k\leq s-1\right\} $ be a sequence of positive numbers, defined by
	$$\beta^s_{s-1}=2^s+2^{s-1}, \ \ \ \beta^s_{s-2}=2^s+2^{s-1}+2^{s-2}, ..., \beta^s_{0}=2^s+2^{s-1}+...+2^{0}.$$
	Then the restricted maximal operator $\widetilde{\sigma }^{\ast ,\nabla }_5,$ defined by 
	\begin{equation*}
	\widetilde{\sigma }^{\ast ,\nabla }_5F:=
\underset{s\in \mathbb{N}}{\sup }\underset{0\leq k\leq s-1}{\sup }\frac{
	\left\vert \sigma _{\beta^s_k}F\right\vert}{\log^2{\beta^s_k}},
	\end{equation*}
	is bounded  from the Hardy space ${{H}_{1/2}}(G)$ to the Lebesgue space ${{L}_{1/2}}(G).$
\end{corollary}
\begin{remark}\label{theorem50}
	Let $f\in {{H}_{1/2}}\left( G \right)$ and 
	$\left\{ \alpha^s_{k}:0\leq k\leq s-1\right\}$ and $\left\{ \beta^s_{k}:0\leq k\leq s-1\right\} $ are   sequences of positive numbers, defined in Corollaries \ref{theorem5} and \ref{theorem6}.
	Then  there exist absolute constants $c_1$ and $c_2$ such that 
	\begin{equation*}
	\left\Vert \sigma _{\alpha^s_k}F\right\Vert_{H_{1/2}}\leq c_1 \left\Vert F\right\Vert_{H_{1/2}} \ \ \ \text{
and \ \ \ }
\left\Vert \sigma _{\beta^s_k}F\right\Vert_{H_{1/2}}\leq c_2 \left\Vert F\right\Vert_{H_{1/2}}
\end{equation*}
for any $ 0\leq k\leq s-1 \ \text{and}\  s\in\mathbb{N}.$

We note that these results was proved in \cite{tep8} and follow the facts that 
$$V(\alpha_k^s)\leq 6 \ \ \ \text{and}\ \ \ V(\beta_k^s)\leq 4, \ \text{ for any} \ \ 0\leq k\leq s-1 \ \text{and}\  s\in\mathbb{N}.$$
\end{remark}

\begin{corollary}\label{theorem7}
	Let $f\in {{H}_{1/2}}\left( G \right)$ and 
	$\left\{ A^s_{k}:0\leq k\leq s-1\right\} $ be a sequence of positive numbers, defined by
	$$\alpha^s_{0}=2^s+2^0, \ \ \ \alpha^s_{1}=2^s+2^0+2^1, ..., \alpha^s_{[\sqrt[n]{s}\ ]-1}=2^s+2^0+2^1+...+2^{[\sqrt[n]{s}\ ]-1}, \ \ \ s\in\mathbb{N},$$
	where $[n]$ denotes the integer part of $n.$
	Then the restricted maximal operator 
	$\widetilde{\sigma }^{\ast ,\nabla }_6,$ defined by 
	\begin{equation*}
	\widetilde{\sigma }^{\ast ,\nabla }_6F:=\underset{s\in \mathbb{N}}{\sup }\underset{0\leq k\leq s-1}{\sup }\frac{
		\left\vert \sigma _{\alpha^s_k}F\right\vert}{\log^{2/n}{\alpha^s_k}},
	\end{equation*}
	is bounded  from the Hardy space ${{H}_{1/2}}(G)$ to the Lebesgue space ${{L}_{1/2}}(G).$
\end{corollary}

\begin{corollary}\label{theorem8}
	Let $f\in {{H}_{1/2}}\left( G \right)$ and 
	$\left\{ A^s_{k}:0\leq k\leq s-1\right\} $ be a sequence of positive numbers, defined by
	$$\beta^s_{s-1}=2^s+2^{s-1}, \ \ \ \beta^s_{s-2}=2^s+2^{-1}+2^{s-2}, ..., \beta^s_{{s-[\sqrt[n]{s}\ ]-1}}=2^s+2^{s-{[\sqrt[n]{s}\ ]}-1}+...+2^{0}.$$
	Then the restricted maximal operator $\widetilde{\sigma }^{\ast ,\nabla }_7,$ defined by 
	\begin{equation*}
	\widetilde{\sigma }^{\ast ,\nabla }_7F:=\underset{s\in \mathbb{N}}{\sup }\underset{0\leq k\leq s-1}{\sup }\frac{
		\left\vert \sigma _{\beta^s_k}F\right\vert}{\log^{2/n}{\beta^s_k}},
	\end{equation*}
	is bounded  from the Hardy space ${{H}_{1/2}}(G)$ to the Lebesgue space ${{L}_{1/2}}(G).$
\end{corollary}

\section{Auxiliary Lemmas}

\begin{lemma}[\textbf{Weisz \cite{We3} and  Simon \cite{S}}]\label{lemma0}
A martingale $F=\left( F_{n},\text{ }n\in \mathbb{N}\right) $
is in $H_{p}\left( 0<p\leq 1\right) $ if and only if there exists a sequence 
$\left( a_{k},k\in \mathbb{N}\right) $ of p-atoms and a sequence $\left( \mu
_{k},k\in \mathbb{N}\right) $ of a real numbers, such that for every $n\in 
\mathbb{N}$
\begin{equation}
\qquad \sum_{k=0}^{\infty }\mu _{k}S_{2^{n}}a_{k}=F_{n},\text{ \ \ \ }%
\sum_{k=0}^{\infty }\left\vert \mu _{k}\right\vert ^{p}<\infty .  \label{6}
\end{equation}
Moreover, 
$\left\Vert F\right\Vert _{H_{p}}\backsim \inf \left(
\sum_{k=0}^{\infty }\left\vert \mu _{k}\right\vert ^{p}\right) ^{1/p},$
where the infimum is taken over all decomposition of $F$ of the form (\ref{6}).
\end{lemma}

\begin{lemma}[\textbf{Weisz \protect\cite{We1}}]
\label{lemma1}Suppose that an operator $T$ is $\sigma $-linear and
\end{lemma}
\begin{equation*}
\int\limits_{\overline{I}}\left\vert Ta\right\vert ^{p}d\mu \leq
c_{p}<\infty ,\text{ \ \ }\left( 0<p\leq 1\right)
\end{equation*}%
for every $p$-atom $a$, where  $I$ denote
the support of the atom. If $T$  is bounded from $L_{\infty \text{ 
}}$  to $L_{\infty },$  then 
\begin{equation*}
\left\Vert TF\right\Vert _{p}\leq c_{p}\left\Vert F\right\Vert _{H_{p}}.
\end{equation*}

\begin{lemma}[see e.g. \protect\cite{G-E-S}, \cite{sws}]\label{lemma2}
Let $t,n\in \mathbb{N}.$ Then 
\begin{equation*}
K_{2^{n}}\left( x\right) =\left\{ 
\begin{array}{c}
\text{ }2^{t-1},\text{\ if \ \ }x\in I_{n}\left( e_{t}\right) ,\text{ }n>t,%
\text{\ }x\in I_{t}\backslash I_{t+1}, \\ 
\left( 2^{n}+1\right) /2,\text{\ if \ \ }x\in I_{n}, \\ 
0,\text{\ otherwise.\ }%
\end{array}%
\right.
\end{equation*}
\end{lemma}

\begin{lemma}[see e.g. \protect\cite{PTWbook}, \cite{tep8}]\label{lemma5}
	Let 
	$$n=\sum_{i=1}^{s}\sum_{k=l_{i}}^{t_{i}}2^{k}, \ \text{ where } \
	t_{1}\geq l_{1}>l_{1}-2\geq t_{2}\geq l_{2}>l_{2}-2>...>t_{s}\geq l_{s}\geq
	0.$$ 
	Then 
	\begin{eqnarray*}
		\left\vert nK_{n}\right\vert &\leq& c\sum_{A=1}^{s}\left( 2^{l_{A}}
		K_{2^{l_{A}}} +2^{t_{A}} K_{2^{t_{A}}}
		+2^{l_A}\sum_{k=l_{A}}^{t_{A}}D_{2^{k}}\right).
	\end{eqnarray*}
\end{lemma}

\begin{lemma}[see e.g. \protect\cite{PTWbook}, \cite{tep8}]\label{corollary1}
Let 
	$$n=\sum_{i=1}^{s}\sum_{k=l_{i}}^{t_{i}}2^{k}, \ \text{where } \ t_{1}\geq l_{1}>l_{1}-2\geq t_{2}\geq l_{2}>l_{2}-2>...>t_{s}\geq l_{s}\geq0.$$ 
	Then, for any $i=1,2,...,s,$	
\begin{equation*}
n\left\vert K_{n}\left( x\right) \right\vert \geq 2^{2t_{i}-5},\text{ \ \
	for \ \ }x\in E_{t_{i}}:=I_{t_{i}+1}\left( e_{t_{i}+1}+e_{t_{i}+2}\right)
\end{equation*}
and
\begin{equation*}
n\left\vert K_{n}\left( x\right) \right\vert \geq 2^{2l_{i}-5},\text{ \ 
	for \  }x\in E_{l_{i}}:=I_{l_{i}+1}\left( e_{l_{i}-1}+e_{l_{i}}\right) , 
\end{equation*}
	where 
	$I_{1}\left( e_{-1}+e_{0}\right) =I_{2}\left( e_{0}+e_{1}\right) .$
\end{lemma}

\section{Proof of the Theorem \ref{theorem1}}

\begin{proof}

Since $\sigma _{n}$ is bounded from $L_{\infty }$ to $%
L_{\infty },$
by Lemma \ref{lemma1}, the proof of theorem \ref{theorem1}
will be complete, if we prove that
\begin{equation*}
\int_{\overline{I_{M}}}\left(\sup_{k\in \mathbb{N}}\left\vert \frac{\sigma _{n_{s_k}}a\left( x\right)}{A^2_{\left\vert n_{s_k}\right\vert}}
\right\vert \right)^{1/2}d\mu \left( x\right) \leq
c<\infty ,
\end{equation*}%
for every 1/2-atom $a.$ We may assume that $a$ is an arbitrary $1/2$-atom,
with support$\ I,$ $\mu \left( I\right) =2^{-M}$ and $I=I_{M}.$ It is easy
to see that 
$\sigma _{n}\left( a\right) =0, \  \text{ when } \  n< 2^{M}.$
Therefore, we can suppose that $n_{s_k}\geq 2^{M}.$
Let $x\in I_{M}$ and $2^{s}\leq n_{s_k}< 2^{s+1}$ for some $n_{s_k}\geq 2^M.$ Since 
$\left\Vert a\right\Vert _{\infty }\leq 2^{2M}$
and 
$\left\vert s_k\right\vert=s,$
by using Lemma \ref{lemma5} we obtain that 
\begin{eqnarray}\label{01}
\ \ \frac{\left\vert \sigma _{n_{s_k}}a\left( x\right) \right\vert}{A^2_{\left\vert n_{s_k}\right\vert}}
&\leq& \frac{2^{2M}}{\left(A_s\right)^2}\int_{I_{M}}\left\vert
K_{n_{s_k}}\left( x+t\right) \right\vert d\mu \left( t\right) \\ \notag
&\leq& \frac{2^{2M}}{\left(A_s\right)^2n_{s_k}}\sum_{A=1}^{r_{s_k}}
\int_{I_{M}}2^{l_{A}^{s_k}}\left\vert K_{2^{l_A^{s_k}}}\left( x+t\right)
\right\vert d\mu \left( t\right)
\\ \notag
&+& \frac{2^{2M}}{\left(A_s\right)^2n_{s_k}}\sum_{A=1}^{r_{s_k}}\int_{I_{M}}2^{t_{A}^{s_k}}\left\vert
K_{2^{t_A^{s_k}}}\left( x+t\right) \right\vert d\mu \left( t\right) 
\\ \notag
&+&\frac{2^{2M}}{\left(A_s\right)^2n_{s_k}}\sum_{A=1}^{r_{s_k}}\int_{I_{M}}2^{l_{A}^{s_k}}
\sum_{k=l_{A}^{s_k}}^{\infty}D_{2^{k}}\left( x+t\right) d\mu \left( t\right)
\\ \notag
&\leq& \frac{2^M}{\left(A_s\right)^22^s}\left(2^M
\sum_{A=1}^{r_s^1}\int_{I_{M}}2^{l_{A}^s}\left\vert K_{2^{l^s_{A}}}\left( x+t\right)
\right\vert d\mu \left( t\right) \right)
\\ \notag
&+& \frac{2^M}{\left(A_s\right)^22^s}\left(2^M\sum_{A=1}^{r^2_{s}}\int_{I_{M}}2^{t^s_{A}}\left\vert
K_{2^{t^s_{A}}}\left( x+t\right) \right\vert d\mu \left( t\right) \right) 
\\ \notag
&+&\frac{2^M}{\left(A_s\right)^22^s}\left(2^M\sum_{A=1}^{r_s^1}\int_{I_{M}}2^{l^s_{j}}
\sum_{k=l^s_j}^{\infty}D_{2^{k}}\left( x+t\right) d\mu \left( t\right) \right).
\end{eqnarray}

If we denote by
\begin{eqnarray*}
	II_{\alpha^s _{A}}^{1}\left( x\right) &:=&2^{M}\int_{I_{M}}2^{\alpha^s
		_{A}}\left\vert K_{2^{\alpha^s _{A}}}\left( x+t\right) \right\vert d\mu \left(
	t\right) ,\text{\ } \alpha=l, \ \text{or} \ \alpha=t\\
	II_{l^s_{A}}^{2}\left( x\right)
	&:=&2^{M}\int_{I_{M}}2^{l^s_{A}}\sum_{k=l^s_{A}}^{\infty}D_{2^{k}}\left(
	x+t\right) d\mu \left( t\right),
\end{eqnarray*}
from \eqref{01} we can conclude that
\begin{eqnarray*}
&&\sup_{2^s\leq n_{s_k}<2^{s+1}}\frac{\left\vert \sigma _{n_{s_k}}a\left( x\right) \right\vert}{A^2_{\left\vert n_{s_k}\right\vert}} \leq\frac{2^M}{\left(A_s\right)^22^s}\left(\sum_{A=1}^{r_s^1} II_{^{l^s_{A}}}^{1}\left(
	x\right) 
	+\sum_{A=1}^{r_s^2}II_{^{t^s_{A}}}^{1}\left( x\right) 
	+\sum_{A=1}^{r_s^1}II_{l_{A}^s}^{2}\left( x\right)
	\right).
\end{eqnarray*}
Hence,

\begin{eqnarray}\label{00}
&&\int_{\overline{I_{M}}}\left(	\sup_{2^s\leq n_{s_k}<2^{s+1}}
\frac{\left\vert \sigma _{n_{s_k}}a\left( x\right) \right\vert}{A^2_{\left\vert n_{s_k}\right\vert}} \right)^{1/2}d\mu(x)\\ \notag
&\leq &\frac{2^{M/2}}{2^{s/2}A_s}\left( \sum_{A=1}^{r^1_s}\int_{\overline{I_{M}}}\left\vert II_{l^s_{A}}^{1}\left( x\right) \right\vert ^{1/2}d\mu \left(	x\right) +\sum_{A=1}^{r^2_s}\int_{\overline{I_{M}}}\left\vert II_{t^s_{A}}^{1}\left( x\right)\right\vert ^{1/2}d\mu \left( x\right) \right.\\ \notag
&&\left.+ \sum_{A=1}^{r^1_s}\int_{\overline{I_{M}}}\left\vert
II_{l^s_{A}}^{2}\left( x\right) \right\vert ^{1/2}d\mu \left( x\right) \right)
\end{eqnarray}

Since 
$$\sup_{s\in \mathbb{N}}r_s\leq A_s \ \ \text{and} \ \ \sup_{s\in \mathbb{N}}r_s^2\leq A_s,$$
we obtain that Theorem \ref{theorem1} is proved if we can prove that 

\begin{equation}\label{11.0}
\int_{%
	\overline{I_{M}}}\left\vert II_{l^s_{A}}^{2}\left( x\right) \right\vert
^{1/2}d\mu \left( x\right) \leq c<\infty , \ \ \ A=1,...,r_s^1 
\end{equation}
and
\begin{equation}\label{11.1}
\int_{\overline{I_{M}}}\left\vert II_{\alpha^s _{A}}^{1}\left( x\right)
\right\vert ^{1/2}d\mu \left( x\right) \leq c<\infty ,
\end{equation}
for all 
$$\alpha^s _{A}=l^s_{A}, \ \ \ A=1,...,r_s^1 \ \ \  and  \ \ \ \alpha^s _{A}=t^s_{A},  \ \ \ A=1,...,r_s^2.$$ 
Indeed, if \eqref{11.0} and \eqref{11.1} hold, from \eqref{00}   we get that

\begin{eqnarray*}
\int_{\overline{I_{M}}}\left(	\sup_{ n_{s_k}\geq2^{M}}\frac{\left\vert \sigma _{n_{s_k}}a\left( x\right) \right\vert}{A^2_{\left\vert n_{s_k}\right\vert}} \right)^{1/2}d\mu(x)
	&\leq&\sum_{s=M}^{\infty}\int_{\overline{I_{M}}}\left(\sup_{2^s\leq n_{s_k}<2^{s+1}}\frac{\left\vert \sigma _{n_{s_k}}a\left( x\right) \right\vert}{A^2_{\left\vert n_{s_k}\right\vert}}\right)^{1/2} \\
	&\leq& \sum_{s=M}^{\infty}\frac{c2^{M/2}}{2^{s/2}}\frac{1}{A_s}\left(2r_s^1+r_s^2\right)\\
	&\leq& \sum_{s=M}^{\infty}\frac{c2^{M/2}}{2^{s/2}}<C<\infty.
\end{eqnarray*}
It remains to prove  \eqref{11.0} and \eqref{11.1}.
Let $t\in I_{M} \  \text{ and } \  x\in I_{l+1}\left( e_{k}+e_{l}\right).$
If
$ 0\leq
k<l<\alpha^s_{A}\leq M \ \ \ \text{or } \  \ \ 0\leq k<l< M< \alpha^s_{A},$ then $x+t\in
I_{l+1}\left( e_{k}+e_{l}\right) $ and if we apply Lemma \ref{lemma2} we obtain that 
\begin{equation}\label{10a}
K_{2^{\alpha^s _{A}}}\left( x+t\right) =0\text{ \ and \ \ }II_{\alpha^s
_{A}}^{1}\left( x\right) =0.  
\end{equation}

Let $0\leq k<\alpha^s _{A}\leq l< M.$
Then $x+t\in I_{l+1}\left( e_{k}+e_{l}\right) $  and if
we use Lemma \ref{lemma2} we get that 
$$2^{\alpha^s _{A}}\left\vert K_{2^{\alpha^s _{A}}}\left( x+t\right) \right\vert\leq 2^{\alpha^s _{A}+k}$$
so that
\begin{equation}\label{10b}
II_{\alpha^s _{A}}^{1}\left(
x\right) \leq 2^{\alpha^s _{A}+k}.  
\end{equation}

Analogously to (\ref{10b}) we can prove that if $0\leq \alpha^s _{A}\leq
k<l< M$, then
\begin{equation*}
2^{\alpha^s _{A}}\left\vert K_{2^{\alpha^s _{A}}}\left( x+t\right) \right\vert
\leq 2^{2\alpha^s _{A}},\text{\ }t\in I_{M},\text{\ }x\in I_{l+1}\left(
e_{k}+e_{l}\right) 
\end{equation*}
so that
\begin{equation}\label{10c}
II_{\alpha^s _{A}}^{1}\left( x\right) \leq
2^{2\alpha^s _{A}},\text{\ }t\in I_{M},\text{\ }x\in I_{l+1}\left(
e_{k}+e_{l}\right) ,  
\end{equation}

Let $t\in I_{M} \  \text{ and }  \   x\in I_M(e_k).$ 
If
$0\leq k<\alpha^s _{A}\leq M$ or $0\leq k<M\leq \alpha^s _{A}.$
Since
$x+t \in x\in I_M(e_k) $ and if we apply Lemma \ref{lemma2}, then we find
that 
$$2^{\alpha^s _{A}}\left\vert K_{2^{\alpha^s _{A}}}\left( x+t\right) \right\vert
\leq 2^{\alpha^s _{A}+k}$$
so that
\begin{equation}
II_{\alpha^s _{A}}^{1}\left(
x\right) \leq 2^{\alpha^s _{A}+k}.  \label{10bb}
\end{equation}
Let
$0\leq \alpha^s _{A}\leq k<M.$
Since
$x+t \in x\in I_M(e_k) $ and if we apply Lemma \ref{lemma2} we obtain
that 
$$2^{\alpha^s _{A}}\left\vert K_{2^{\alpha^s _{A}}}\left( x+t\right) \right\vert
\leq 2^{2\alpha^s _{A}}$$
and
\begin{equation}\label{10bbb}
II_{\alpha^s _{A}}^{1}\left(
x\right) \leq 2^{2\alpha^s _{A}}.  
\end{equation}

Let $0\leq \alpha^s _{A}< M, A=1,...,s.$ By combining (\ref{1}) with (\ref{10a})-(\ref{10bbb})  we get that

\begin{eqnarray*}
&&\int_{\overline{I_{M}}}\left\vert II_{\alpha^s _{A}}^{1}\left( x\right)
\right\vert ^{1/2}d\mu \left( x\right) \\
&=&\overset{M-2}{\underset{k=0}{\sum }}\overset{M-1}{\underset{l=k+1}{\sum }}
\int_{I_{l+1}\left( e_{k}+e_{l}\right) }\left\vert II_{\alpha^s
_{A}}^{1}\left( x\right) \right\vert ^{1/2}d\mu \left( x\right) +\overset{M-1
}{\underset{k=0}{\sum }}\int_{I_{M}\left( e_{k}\right) }\left\vert
II_{\alpha^s _{A}}^{1}\left( x\right) \right\vert ^{1/2}d\mu \left( x\right) \\
&\leq&c\overset{\alpha^s _{A}-1}{\underset{k=0}{\sum }}\overset{M-1}{\underset{%
l=\alpha^s _{A}}{\sum }}\int_{I_{l+1}\left( e_{k}+e_{l}\right) }2^{\left(
\alpha^s _{A}+k\right) /2}d\mu \left( x\right) +c\overset{M-2}{\underset{%
k=\alpha^s _{A}}{\sum }}\overset{M-1}{\underset{l=k+1}{\sum }}%
\int_{I_{l+1}\left( e_{k}+e_{l}\right) }2^{\alpha^s _{A}}d\mu \left( x\right) \\
&+&c\overset{\alpha^s _{A}-1}{\underset{k=0}{\sum }}\int_{I_{M}\left(
e_{k}\right) }2^{\left( \alpha^s _{A}+k\right) /2}d\mu \left( x\right) +c%
\overset{M-1}{\underset{k=\alpha^s _{A}}{\sum }}\int_{I_{M}\left( e_{k}\right)
}2^{\alpha^s _{A}}d\mu \left( x\right) \\
&\leq& c\overset{\alpha^s _{A}-1}{\underset{k=0}{\sum }}\overset{M-1}{\underset{%
l=\alpha^s _{A}+1}{\sum }}\frac{2^{\left( \alpha^s _{A}+k\right) /2}}{2^{l}}+c%
\overset{M-2}{\underset{k=\alpha^s _{A}}{\sum }}\overset{M-1}{\underset{l=k+1}{%
\sum }}\frac{2^{\alpha^s _{A}}}{2^{l}}\\
&+&c\overset{\alpha^s _{A}-1}{\underset{k=0}{
\sum }}\frac{2^{\left( \alpha^s _{A}+k\right) /2}}{2^{M}}+c\overset{M-1}{%
\underset{k=\alpha^s _{A}}{\sum }}\frac{2^{\alpha^s _{A}}}{2^{M}}
\leq C<\infty.
\end{eqnarray*}

Analogously we can prove that (\ref{11.1}) holds also for the case $\alpha^s _{A}\geq M.$ Hence, \eqref{11.1} holds and it remains to prove \eqref{11.0}.

Now, prove boundedness of $II_{l^s_{A}}^{2}$. Let 
$t\in I_{M} \  \text{ and } \  x\in
I_{i}\backslash I_{i+1}.$ \
If $i\leq l^s_{A}-1,$ since $x+t\in I_{i}\backslash
I_{i+1},$ by using (\ref{1dn}) we have that 
\begin{equation}\label{13a}
II_{l^s_{A}}^{2}\left( x\right) =0.  
\end{equation}

If $l^s_{A}\leq i< M,$  by using (\ref{1dn}) we obtain that 
\begin{equation}\label{13b}
II_{l^s_{A}}^{2}\left( x\right) \leq
2^{M}\int_{I_{M}}2^{l^s_{A}}\sum_{k=l^s_{A}}^{i}D_{2^{k}}\left( x+t\right) d\mu
\left( t\right) \leq c2^{l^s_{A}+i}.  
\end{equation}

Let $0\leq l^s_{A}< M.$ By combining (\ref{1}), (\ref{13a}) and (\ref
{13b}) we get that 
\begin{eqnarray}\label{star1}
&&\int_{\overline{I_{M}}}\left\vert II_{l^s_{A}}^{2}\left( x\right) \right\vert
^{1/2}d\mu \left( x\right)\\ \notag
 &=&\left(
\sum_{i=0}^{l^s_{A}-1}+\sum_{i=l^s_{A}+1}^{M-1}\right)
\int_{I_{i}\backslash I_{i+1}}\left\vert II_{l^s_{A}}^{2}\left( x\right)
\right\vert ^{1/2}d\mu \left( x\right)\\ \notag
&\leq& c\sum_{i=l^s_{A}}^{M-1}\int_{I_{i}\backslash I_{i+1}}2^{\left(
l^s_{A}+i\right) /2}d\mu \left( x\right)\\ \notag
&\leq& c\sum_{i=l^s_{A}}^{M-1}2^{\left( l^s_{A}+i\right) /2}\frac{1}{2^{i}}\leq
C<\infty .
\end{eqnarray}

If  $M\leq l^s_{A},$ then $i<M\leq l^s_{A}$ and apply \eqref{13a} we get that
\begin{eqnarray}\label{star2}
\int_{\overline{I_{M}}}\left\vert II_{l^s_{A}}^{2}\left( x\right) \right\vert
	^{1/2}d\mu \left( x\right)=0,
\end{eqnarray}
and also \eqref{11.0} is proved by just combining \eqref{star1} and \eqref{star2} so part a) is complete and we turn to the proof of b).

Under condition (\ref{cond2}), there exists an increasing sequence $\left\{ \alpha _{k}:\text{ }k\geq
0\right\} \subset \left\{ n_{k}\text{ }:k\geq 0\right\} $ of positive
integers, such that 
\begin{equation*}
\lim_{k\to\infty}\frac{\left(\vert A_{\left\vert\alpha _{k}\right\vert}\vert\right)^{1/2}}{\varphi^{1/4}_{\left\vert\alpha _{k}\right\vert}}=\infty
\end{equation*}
and
\begin{equation}\label{2aaa}
\sum_{k=1}^{\infty }\frac{\varphi^{1/4}_{\left\vert\alpha _{k}\right\vert}}{\left(\vert A_{\left\vert\alpha _{k}\right\vert}\vert\right)^{1/2}} \leq c<\infty .  
\end{equation}

Let 
\begin{equation*}
F_{A}:=\sum_{\left\{ k;\text{ }\left\vert \alpha _{k}\right\vert <A\right\}
}\lambda _{k}a_{k},
\end{equation*}
where
\begin{equation*}
\lambda _{k}:=\frac{\varphi^{1/2}_{\left\vert\alpha _{k}\right\vert}}{\vert A_{\left\vert\alpha _{k}\right\vert}\vert} \ \ \   \text{
and} \  \ \ 
a_{k}:=2^{\left\vert \alpha _{k}\right\vert }\left(
D_{2^{\left\vert \alpha _{k}\right\vert +1}}-D_{2^{\left\vert \alpha
		_{k}\right\vert }}\right) .
\end{equation*}

Since \ \ 
$$
\text{supp}(a_{k})=I_{\left\vert \alpha _{k}\right\vert },
\ \ \ 
\left\Vert
a_{k}\right\Vert _{\infty }\leq 2^{2\left\vert \alpha _{k}\right\vert }=\mu (%
\text{supp }a_{k})^{-2} \ \
$$ 
and
\begin{equation*}
S_{2^{A}}a_{k}=\left\{ 
\begin{array}{ll}
a_{k} & \left\vert \alpha _{k}\right\vert <A, \\ 
0\, & \left\vert \alpha _{k}\right\vert \geq A,%
\end{array}%
\right.  
\end{equation*}
if we apply Lemma \ref{lemma0} and (\ref{2aaa}) we can conclude that $F=\left(
F_{1},F_{2},...\right) \in H_{1/2}.$

It is easy to prove that
\begin{equation}\label{5aa}
\widehat{F}(j)  
=\left\{ 
\begin{array}{ll}
\left(2^{\left\vert \alpha _{k}\right\vert } {\varphi^{1/2}_{\left\vert\alpha _{k}\right\vert}} \right)
/\left\vert A_{\left\vert\alpha _{k}\right\vert}\right\vert , & \text{\thinspace \thinspace }j\in \left\{
2^{\left\vert \alpha _{k}\right\vert },...,2^{_{\left\vert \alpha
_{k}\right\vert +1}}-1\right\} ,\text{ }k=0,1,... \\ 
0\,, & \text{\thinspace }j\notin \bigcup\limits_{k=0}^{\infty }\left\{
2^{_{\left\vert \alpha _{k}\right\vert }},...,2^{_{\left\vert \alpha
_{k}\right\vert +1}}-1\right\} .\text{ }%
\end{array}%
\right.
\end{equation}

Let $2^{\left\vert \alpha _{k}\right\vert }<j<\alpha _{k}.$ By using (\ref%
{5aa}) we get that%
\begin{eqnarray}\label{sn} \ \ \ \ \
S_{j}F=S_{2^{\left\vert \alpha _{k}\right\vert }}F+\sum_{v=2^{^{\left\vert
\alpha _{k}\right\vert }}}^{j-1}\widehat{F}(v)w_{v} =S_{2^{\left\vert \alpha
_{k}\right\vert }}F+\frac{\left( D_{j}-D_{2^{\left\vert \alpha
_{k}\right\vert }}\right) 2^{\left\vert \alpha _{k}\right\vert }}{\left\vert A_{{\left\vert\alpha _{k}\right\vert}}\right\vert}  
\end{eqnarray}
Let \
$\begin{matrix}
2^{\left\vert\alpha _{k}\right\vert}\le {{\alpha}_{{{s}_{n}}}}\le {2^{{\left\vert\alpha _{k}\right\vert}+1}}.
\end{matrix}$ \
Then, by using \eqref{sn} we find that
\begin{eqnarray}\label{7aaa} 
&&\sigma_{\alpha_{s_n}}F\\ \notag
&=&\frac{1}{\alpha_{s_n}}\sum_{j=1}^{2^{\left\vert \alpha_{k}\right\vert }}S_{j}F+\frac{1}{\alpha _{s_n}}
\sum_{j=2^{\left\vert \alpha_{k}\right\vert }+1}^{\alpha _{s_n}}S_{j}F \\ \notag
&=&\frac{\sigma_{2^{\left\vert \alpha _{k}\right\vert }}F}{\alpha _{s_n}}+\frac{\left( \alpha _{k}-2^{\left\vert
\alpha _{k}\right\vert }\right) S_{2^{\left\vert \alpha _{k}\right\vert }}F}{ \alpha _{s_n}}+\frac{2^{\left\vert \alpha
_{k}\right\vert } {\varphi^{1/2}_{\left\vert\alpha _{k}\right\vert}}}{\left\vert A_{{\left\vert\alpha _{k}\right\vert}}\right\vert \alpha _{s_n}}\sum_{j=2^{_{\left\vert
\alpha _{k}\right\vert }}+1}^{\alpha_{s_n}}\left( D_j-D_{2^{\left\vert
\alpha _{k}\right\vert }}\right)\\ \notag
&:=&III_{1}+III_{2}+III_{3}.
\end{eqnarray}

Since
$$D_{j+2^{m}}=D_{2^{m}}+w_{_{2^{m}}}D_{j},\text{ \  when \ }j<2^{m}  $$
we obtain that

\begin{eqnarray}\label{9aaa} \ \ \ \ \ 
\left\vert III_{3}\right\vert &=&\frac{2^{\left\vert \alpha _{k}\right\vert
}{\varphi^{1/2}_{\left\vert\alpha _{k}\right\vert}}}{\left\vert A_{\left\vert\alpha _{k}\right\vert}\right\vert \alpha _{s_n}}\left\vert \sum_{j=1}^{\alpha
_{s_n}-2^{_{\left\vert \alpha_{k}\right\vert }}}\left( D_{j+2^{_{\left\vert
\alpha _{k}\right\vert }}}-D_{2^{\left\vert \alpha _{k}\right\vert }}\right)
\right\vert \\ \notag
&=&\frac{2^{\left\vert \alpha _{k}\right\vert }{\varphi^{1/2}_{\left\vert\alpha _{k}\right\vert}}}{\left\vert A_{\left\vert\alpha _{k}\right\vert}\right\vert \alpha _{s_n}}\left\vert \sum_{j=1}^{\alpha _{s_n}-2^{\left\vert \alpha_{k}\right\vert }}D_j\right\vert \\ \notag
&=&\frac{2^{\left\vert \alpha_{k}\right\vert } {\varphi^{1/2}_{\left\vert\alpha _{k}\right\vert}}}{\left\vert A_{\left\vert\alpha _{k}\right\vert}\right\vert \alpha_{s_{n}}}\left( \alpha
_{s_n}-2^{\left\vert \alpha_{k}\right\vert }\right) \left\vert K_{\alpha
_{s_n}-2^{\left\vert \alpha_{k}\right\vert }}\right\vert
\\ \notag
&\geq&\frac{{\varphi^{1/2}_{\left\vert\alpha _{k}\right\vert}}}{2\left\vert A_{\left\vert\alpha _{k}\right\vert}\right\vert }\left( \alpha_{s_n}-2^{\left\vert \alpha _{k}\right\vert }\right) \left\vert K_{\alpha_{s_{n}}-2^{\left\vert \alpha _{k}\right\vert }}\right\vert.
\end{eqnarray}
By combining well-known estimates (see \cite{PTWbook})
	\begin{equation*}
\left\Vert S_{2^k}F\right\Vert_{H_{1/2}}\leq c_1 \left\Vert F\right\Vert_{H_{1/2}} \ \ \ \text{and } \ \ \ \left\Vert \sigma_{2^k}F\right\Vert_{H_{1/2}}\leq c_2 \left\Vert F\right\Vert_{H_{1/2}}, \ \ \ \ k\in \mathbb{N},
\end{equation*}
we obtain that
$$\left\Vert III_{1}\right\Vert
_{1/2}\leq C \ \ \ \text{and} \ \ \ \left\Vert III_{2}\right\Vert _{1/2}\leq C.$$ 

Let \
$\begin{matrix}
2^{\left\vert\alpha _{k}\right\vert}\le {{\alpha}_{{{s}_{1}}}}\le {{\alpha}_{{{s}_{2}}}}\le ...\le {{\alpha}_{{{s}_{r}}}}\le {2^{{\left\vert\alpha _{k}\right\vert}+1}}
\end{matrix}$ \
be natural numbers which generates the set 
$$
{{A}_{\left\vert\alpha _{k}\right\vert}}={\left\{ l_{1}^{\left\vert\alpha _{k}\right\vert},l_{2}^{\left\vert\alpha _{k}\right\vert},...,l_{r^1_{\left\vert\alpha _{k}\right\vert}}^{\left\vert\alpha _{k}\right\vert} \right\}}
\bigcup{\left\{ t_{1}^{\left\vert\alpha _{k}\right\vert},t_{2}^{\left\vert\alpha _{k}\right\vert},...,t_{{{r}_{\left\vert\alpha _{k}\right\vert}^2}}^{\left\vert\alpha _{k}\right\vert} \right\}}
$$
and choose number $\alpha_{s_{n}}=\sum_{i=1}^{r_{n}}\sum_{k=l_{i}^{n}}^{t_{i}^{n}}2^{k},$ where
\begin{eqnarray*}
t_{1}^{\left\vert\alpha _{k}\right\vert}\geq l_{1}^{\left\vert\alpha _{k}\right\vert}>l_{1}^{\left\vert\alpha _{k}\right\vert}-2\geq t_{2}^{\left\vert\alpha _{k}\right\vert}\geq
l_{2}^{\left\vert\alpha _{k}\right\vert}>l_{2}^{\left\vert\alpha _{k}\right\vert}-2\geq ...\geq t_{\left\vert\alpha _{k}\right\vert}^{\left\vert\alpha _{k}\right\vert}\geq l_{\left\vert\alpha _{k}\right\vert}^{\left\vert\alpha _{k}\right\vert}\geq 0,
\end{eqnarray*}
for some $1\leq n\leq r,$ such that 
$l^{\left\vert\alpha _{k}\right\vert}_u=l_i,  \  \text{for some} \ 1\leq u\leq r^1_{\left\vert\alpha _{k}\right\vert},  1\leq i\leq r_{\left\vert\alpha _{k}\right\vert}^1. $

Since
$\mu \left\{E_{l_{i}}\right\} \geq 1/2^{l_i+1},$
by using Lemma \ref{corollary1} we get that 
\begin{eqnarray}\label{low1}
\int_{E_{l_{i}}}\left\vert\underset{k\in \mathbb{N}}{\sup }
\frac{\left\vert \sigma _{\alpha _{s_k}}F\right\vert}{\varphi_{\vert \alpha _{s_k}\vert}}\right\vert^{1/2}d\mu 
&\geq&	\int_{E_{l_{i}}}\left\vert \frac{\sigma _{\alpha _{s_n}}F(x)}{{{\varphi_{\left\vert\alpha _{s_n}\right\vert}}}}\right\vert^{1/2} d\mu \\ \notag
&\geq& \frac{2^{\left(2l_{i}-6\right)/2}}{\sqrt{2}\left(\left\vert A_{ {\left\vert\alpha _{k}\right\vert}}\right\vert \right)^{1/2}{\varphi^{1/4}_{\left\vert\alpha _{k}\right\vert}}}\frac{1}{2^{l_i+1}}\\ \notag
&\geq&\frac{1}{2^5\left(\left\vert A_{{\left\vert\alpha _{k}\right\vert}}\right\vert\right)^{1/2}{\varphi^{1/4}_{\left\vert\alpha _{k}\right\vert}}}.
\end{eqnarray}

On the other hand, we can also choose number $\alpha_{s_{n}},$ for some $1\leq n\leq r,$ such that 
$t^{\left\vert\alpha _{k}\right\vert}_u=t_i,  \  \text{for some} \ 1\leq u\leq r^2_{\left\vert\alpha _{k}\right\vert}, \ \  1\leq i\leq r_{\left\vert\alpha _{k}\right\vert}^2. $
According the fact that 
$\mu \left\{
E_{t_{i}}\right\} \geq 1/2^{t_{i}+3},$
by using again Lemma \ref{corollary1} for some $\alpha _{k}$ and $1\leq i\leq r_s^2$ we also get that

\begin{eqnarray}\label{low2} \ \ \  \ \ \
	\int_{E_{l_{i}}}\left\vert\underset{k\in \mathbb{N}}{\sup }
	\frac{\left\vert \sigma _{\alpha _{s_k}}F\right\vert}{\varphi_{\vert \alpha _{s_k}\vert}}\right\vert^{1/2}d\mu 
	&\geq&	\int_{E_{t_{i}}}\left\vert \frac{\sigma _{\alpha _{s_n}}F(x)}{{{\varphi_{\left\vert\alpha _{s_n}\right\vert}}}}\right\vert^{1/2} d\mu \\ \notag
	&\geq& \frac{1}{\sqrt{2}\left(\left\vert A_{{\left\vert\alpha _{k}\right\vert}}\right\vert\right)^{1/2}{\varphi^{1/4}_{\left\vert\alpha _{k}\right\vert}} }2^{\left(2t_{i}-6\right)/2}\frac{1}{2^{t_i+3}}\\ \notag
	&\geq&\frac{1}{2^7\left(\left\vert A_{{\left\vert\alpha _{k}\right\vert}}\right\vert\right)^{1/2}{\varphi^{1/4}_{\left\vert\alpha _{k}\right\vert}}}.
\end{eqnarray}
By combining (\ref{7aaa})-(\ref{low2}) and Lemma \ref{corollary1} for sufficiently big $\alpha_k$ we obtain that 
\begin{eqnarray*}
&&\int_{G}\left\vert    \underset{k\in \mathbb{N}}{\sup }
\frac{\left\vert \sigma _{\alpha _{s_k}}F\right\vert}{\varphi_{\vert \alpha _{s_k}\vert}}       \right\vert^{1/2} d\mu \\
&\geq&\left\Vert III_{3}\right\Vert
_{1/2}^{1/2}-\left\Vert III_{2}\right\Vert _{1/2}^{1/2}-\left\Vert
III_{1}\right\Vert _{1/2}^{1/2}\\
&\geq& \underset{i=1}{\overset{r^1_{\left\vert\alpha _{k}\right\vert}-1}{\sum }}\int_{E_{l_i}}\left\vert \underset{k\in \mathbb{N}}{\sup }
\frac{\left\vert \sigma _{\alpha _{s_k}}F\right\vert}{\varphi_{\vert \alpha _{s_k}\vert}}\right\vert^{1/2} d\mu +\underset{i=1}{\overset{r^2_{\left\vert\alpha _{k}\right\vert}-1}{\sum }}\int_{E_{t_i}}\left\vert \underset{k\in \mathbb{N}}{\sup }
\frac{\left\vert \sigma _{\alpha _{s_k}}F\right\vert}{\varphi_{\vert \alpha _{s_k}\vert}}\right\vert^{1/2} d\mu -2C\\
&\geq&\frac{r_{\left\vert\alpha _{k}\right\vert}^1
	+r_{\left\vert\alpha _{k}\right\vert}^2}{2^7\left(\left\vert A_{{\left\vert\alpha _{k}\right\vert}}\right\vert\right)^{1/2}\varphi^{1/4}_{\left\vert\alpha _{k}\right\vert}}-2C
\geq \frac{\left(\left\vert A_{\left\vert\alpha _{k}\right\vert}\right\vert\right)^{1/2}}{2^8\varphi^{1/4}_{\left\vert\alpha _{k}\right\vert}}\to \infty, \ \ \ \text{as} \ \ \ k\to \infty.
\end{eqnarray*}
so also Part b) id proved and 
the proof is complete.
\end{proof}


\begin{thebibliography}{99}
	
\bibitem{BNPT} \textit{D. Baramidze, N. Nadirashvili, L.-E. Persson} and \textit{G. Tephnadze,} Some weak-type inequalities and almost everywhere convergence of Vilenkin-Nörlund means, J. Inequal. Appl. 2023 (to appear).

\bibitem{BaTe1} \textit{D. Baramidze} and \textit{G. Tephnadze,} Restricted maximal operators of Fejér means of Walsh-Fourier series in the space $H_{1/2}$,  Banach J. Math. Anal. (to appear).

\bibitem{Fu} \textit{N. J. Fujii,} A maximal inequality for $H_{1}$
functions on the generalized Walsh-Paley group, Proc. Amer. Math. Soc. 77
(1979), 111-116.

\bibitem{gat1} \textit{G. Gát,} Inverstigations of certain operators with
respect to the Vilenkin sistem, Acta Math. Hung. 61 (1993), 131-149.

\bibitem{GoAMH} \textit{U. Goginava,} Maximal operators of Fejér means of
double Walsh-Fourier series. Acta Math. Hungar. 115 (2007), 333--340.

\bibitem{Goginava} \textit{U. Goginava,} The martingale Hardy type
inequality for Marcinkiewicz-Fejér means of two-dimensional conjugate
Walsh-Fourier series, Acta Math. Sinica 27 (2011), 1949-1958.

\bibitem{GoSzeged} \textit{U. Goginava,} Maximal operators of Fejér-Walsh
means. Acta Sci. Math. (Szeged) 74 (2008), 615--624.

\bibitem{G-E-S} \textit{B. Golubov, A. Efimov} and \textit{V. Skvortsov,} Walsh series and transformations, Dordrecht, Boston, London, 1991. Kluwer
Acad. Publ., 1991.

\bibitem{NTT} \textit{ N. Nadirashvili, G. Tephnadze}  and \textit{  G. Tutberidze,} Almost everywhere and norm convergence of Approximate Identity and Fejér means of trigonometric and Vilenkin systems, Trans. A. Razmadze Math. Inst. 2023 (to appear).

\bibitem{PS} \textit{J. Pál } and \textit{P. Simon,} On a generalization of
the concept of derivate, Acta Math. Hung. 29 (1977), 155-164.


\bibitem{PT} \textit{L. E. Persson} and  \textit{G. Tephnadze,} A sharp boundedness result concerning some maximal operators of Vilenkin-Fejér means, Mediterr. J. Math. 13, 4 (2016) 1841-1853.

\bibitem{PTW2} \textit{L. E. Persson, G. Tephnadze} and  \textit{P. Wall,} On the maximal operators of Vilenkin-Nörlund means, J. Fourier Anal. Appl. 21, 1 (2015), 76-94. 

\bibitem{PTWbook} \textit{L. E. Persson, G. Tephnadze} and \textit{F. Weisz,} Martingale Hardy Spaces and Summability of Vilenkin-Fourier Series, book manuscript, Birkhäuser, 2022, (to appear).

\bibitem{sws} \textit{F. Schipp, W. Wade, P. Simon} and \textit{ J. Pál, } Walsh series, An Introduction to Duadic Harmonic Analysis, Akademiai Kiadó, 
Budapest-Adam Hilger, 1990.

\bibitem{Sc} \textit{F. Schipp,} Certain rearrangements of series in the
Walsh series, Mat. Zametki, 18 (1975), 193-201.

\bibitem{Si1} \textit{P. Simon,} Cesáro summability with respect to
two-parameter Walsh systems, Monatsh. Math. 131, 4 (2000), 321--334.

\bibitem{Si2} \textit{P. Simon,} Investigations with respect to the Vilenkin
system, Ann. Univ. Sci. Budapest. Eötvös Sect. Math. 28 (1985), 87-101.


\bibitem{S} \textit{P. Simon,} A note on the of the Sunouchi operator with
respect to Vilenkin systems, Annales Univ. Sci. Sect. Math. Budapest (2001).

\bibitem{tep1} \textit{G. Tephnadze,} Fejér means of Vilenkin-Fourier
series, Stud. sci. math. Hung. 49, (2012) 79-90.

\bibitem{tep2} \textit{G. Tephnadze,} On the maximal operator of Vilenkin-Fejér means, Turk. J. Math. 37, (2013), 308-318.

\bibitem{tep3} \textit{G. Tephnadze,} On the maximal operators of
Vilenkin-Fejér means on Hardy spaces, Math. Inequal. Appl. 16,
2 (2013),  301-312.


\bibitem{tep8} \textit{G. Tephnadze,} On the convergence of Fejér means of Walsh-Fourier series in the space $H_p$, J. Contemp. Math. Anal. 51, 2 (2016), 90-102.


\bibitem{We2} \textit{F. Weisz}, Cesáro summability of one- and
two-dimensional Walsh-Fourier series, Anal. Math. 22 (1996), 229-242.

\bibitem{We1} \textit{F. Weisz,} Martingale Hardy spaces and their
applications in Fourier Analysis, Springer, Berlin-Heideiberg-New York, 1994.

\bibitem{We3} \textit{F. Weisz,} Summability of multi-dimensional Fourier
series and Hardy space, Kluwer Academic, Dordrecht, 2002.

\bibitem{We4} \textit{F. Weisz,} Weak type inequalities for the Walsh and
bounded Ciesielski systems, Anal. Math. 30,  2 (2004),  147-160.
\end{thebibliography}
\end{document}